\documentclass[10pt, reqno, fleqn, letter]{amsart}
\usepackage[french]{babel}
\usepackage[T1]{fontenc}
\usepackage[latin1]{inputenc}
\usepackage{array}
\usepackage{amsfonts}
\usepackage{amssymb}
\usepackage{amsmath}
\usepackage{amsthm}
\usepackage{mathrsfs}
\usepackage[all]{xy}
\usepackage{fancyhdr}
\usepackage{a4wide}
\usepackage{hyperref}
\usepackage[bitstream-charter]{mathdesign} %In this one, replace \nu with \!\nu
\newtheorem*{thm*}{Théorème}

\newcommand{\R}{\mathbb{R}}
\newcommand{\Z}{\mathbb{Z}}

\newcommand{\Tr}{\mathrm{Tr}}
\newcommand{\Str}{\mathrm{Tr}_s}

\newcommand{\Ind}{\mathrm{Ind}}
\newcommand{\End}{\mathrm{End}}
\newcommand{\ch}{\mathrm{ch}}
\newcommand{\sh}{\mathrm{sh}}

\newcommand{\ad}{\mathrm{ad}}

\renewcommand{\dim}{\mathrm{dim}}

\renewcommand{\det}{\mathrm{det}}

\title{Calcul du cocycle JLO pour l'opérateur de Dirac classique} 
\author{ 
Rudy Rodsphon 
}
\date{31 août 2011}

\begin{document}

\maketitle 

\setlength{\parindent}{0mm}

\section{Introduction}

Le théorème de l'indice, démontré en 1963 par M. Atiyah et I. Singer, est certainement l'un des théorèmes les plus marquants du vingtième siècle. Il crée un pont entre diverses branches des mathématiques et a de nombreuses applications. Par exemple, les théorèmes de Gauss-Bonnet et Riemann-Roch-Hirzebruch sont des cas particuliers de théorème de l'indice. \\

Le théorème de l'indice fait le lien entre topologie et analyse, plus précisément, il montre que l'indice analytique d'un opérateur différentiel elliptique (donc de Fredholm) sur une variété compacte $X$  coïncide avec son indice topologique en K-théorie, qui est une application $\Ind : K_0(X) \longrightarrow \Z $. Cependant, la richesse de la théorie de l'indice ne provient pas seulement du résultat en lui-même, mais également des techniques utilisées pour le démontrer. Comme on peut s'y attendre, les ingrédients de la première preuve étaient essentiellement dûs à des arguments de topologie différentielle (théorie du cobordisme) et de K-théorie, permettant de ramener les hypothèses générales à des cas particuliers où l'on sait calculer directement les indices en question, ce qui rendait sa preuve particulièrement difficile, bien que conceptuelle.  Beaucoup plus surprenant, une preuve du théorème de l'indice utilisant à l'équation de la chaleur a été exposée pour la première fois par M. Atiyah, R. Bott et V.K. Patodi en 1973 [A-B-P], puis sur cette base, E. Getzler a simplifié cette preuve (\cite{Ge}, \cite{BGV}), utilisant des méthodes venant la théorie quantique des champs initiées par Alvarez-Gaumé et E. Witten. Plus récemment, le développement de la géométrie non commutative a également contribué à apporter des idées nouvelles à ce sujet. \\

Rappelons tout d'abord le théorème de l'indice pour les opérateurs de Dirac, utilisant l'équation de la chaleur et le calcul de Getzler développé dans \cite{Ge}, puis dans \cite{BGV} avec divers remaniements. Enonçons le théorème en question. 

\begin{thm*}
Soit $M$ une variété riemannienne fermée orientée de dimension paire $n$, $E$ un module de Clifford au-dessus de $M$ et $D$ un opérateur de Dirac sur $E$. Alors l'indice de $D$ est donné par la formule suivante
$$ \Ind(D) = \frac{1}{(2 i\pi)^{n/2}}\int_M \widehat{A}(M) \ch(E/S) $$
\end{thm*}
\noindent où $\widehat{A}(M)$ est la classe caractéristique donnée par la formule 
$$\widehat{A}(M) = \det^{1/2} \left( \frac{R/2}{\sh(R/2)} \right) \in \Omega^{4\bullet}(M)$$
$R \in \Omega^{2\bullet}(M, End(T^*M))$ étant la courbure de Riemann de $M$.  \\

A défaut de redonner la définition du caractère de Chern tordu $\ch(E/S)$, rappelons que si $M$ est aussi spin, et que l'on a $E = S \otimes W$, où $S$ est le fibré des spineurs et $W$ un fibré vectoriel sur $M$, alors $\ch(E/S)=\ch(W)$. Dans le cas général, une telle décomposition n'est possible que localement, $\ch(E/S)$ permet de contourner ce problème. \\

Rappelons les étapes et idées clés de la preuve. Nous travaillons dans le cadre suivant : soient $M$ une variété riemannienne compacte sans bord, $E$ un module de Clifford au-dessus de $M$, c'est à dire un fibré vectoriel où chaque fibre en $x \in M$ est une représentation de l'algèbre de Clifford $C(T_x^*M)$, et $D$ un opérateur de Dirac sur $E$. $D^2$ étant un Laplacien généralisé, la première étape étant de voir que l'indice de $D$ est donné par la formule suivante  (McKean-Singer)
$$ \Ind(D) = \Str(e^{-tD^2}) = \int_M \Str(p_t(x,x)) \vert dx \vert$$
où $e^{-tD^2}$ est l'opérateur de la chaleur associé au noyau de la chaleur $p_t(x,y) \in \Gamma(E \boxtimes E^*)$, $t \in \R_+$, vérifiant l'équation de la chaleur $(\partial_t + D^2_x)p_t(x,y)  = 0$, pour tout $x,y \in M$. \\

L'idée qui transparaît tout au long de cette preuve est que l'on quantifie des objets géométriques à partir des théories classiques. Les opérateurs de Dirac sont au départ une quantification des connexions : en effet, classiquement, étant donné une connexion $\nabla^E$ sur $E$, un opérateur de Dirac se définit comme la composition 
$$ D : \Gamma(E) \overset{\nabla^E}{\longrightarrow} \Gamma(T^*M \otimes E) \overset{c}{\longrightarrow} \Gamma(E)$$
et inversement, on peut montrer que tout opérateur de Dirac généralisé est associé à une superconnexion. La formule de Lichnerowicz met en relation le carré de l'opérateur de Dirac et le Laplacien associé à la connexion $\nabla^E$, et fait ainsi apparaître $\Str(e^{-tD^2})$ comme une quantification du caractère de Chern $\ch(\nabla^E) = \Str(e^{-F^E})$, où $F^E$ est la courbure de la connexion $\nabla^E$. De même, le fibré de Clifford est une quantification du fibré des formes différentielles sur $M$. \\

Revenons au noyau de la chaleur $p_t$. Fixons $x_0 \in M$, et paramétrisons $M$ en coordonnées normales dans un voisinage de $x_0$ suffisamment petit. Après s'être mis dans une trivialisation locale de $E$ convenable, le noyau de la chaleur $p_t$ peut-être vu comme une fonction $v \mapsto k(t,v)$ sur un voisinage $U$ de $0$ dans $T_{x_0}M$, à valeurs dans 
$$ \End_{x_0}(E) = C(T_{x_0}^*M) \otimes \End_{C(T_{x_0}^*M)}(E) $$
Le symbole $\sigma$ donne un isomorphisme $C(T_{x_0}^*M) \longrightarrow \Lambda^\bullet T_{x_0}M$, $\Lambda^\bullet T_{x_0}M$ étant un module sur l'algèbre de Clifford $C(T_{x_0}^*M)$, et ainsi l'opérateur de la chaleur $\partial_t + D^2$, à coefficients dans $C(T_{x_0}^*M) \otimes \End_{C(T_{x_0}^*M)}(E)$, agit sur $k(t,v)$ alors vu comme élément de $C^{\infty}(U,\Lambda^\bullet T_{x_0}M \otimes \End_{C(T_{x_0}^*M)}(E))$ au travers du symbole $\sigma$.  \\

L'idée qui va nous permettre de conclure va être d'introduire un réechelonnement des fonctions sur $\R_+ \times U$ à valeurs dans $\Lambda^\bullet T_{x_0}M \otimes \End_{C(T_{x_0}^*M)}(E)$ en posant 
$$ r(u,t,v) = u^{n/2}(\delta_u k)(t,v) = \sum_{i=0}^n u^{(n-i)/2} k(ut,u^{1/2}v)_{[i]} \, , \, (t,v) \in \R_+ \times U $$
de sorte à avoir 
$$\underset{u \to 0}{\lim} \, r(u,t=1,v=0) = \underset{u \to 0}{\lim} \sum_{i=0}^n u^{(n-i)/2} k_u(x_0,x_0)_{[i]} $$
L'indice $[i]$ signifie que l'on extrait la partie de degré $i$ de la forme differentielle à laquelle la notation s'applique. Lorsque $u \to 0$, il ne reste que la composante de degré $n$, qui est en fait la seule qui nous intéresse puisque la supertrace annule les composantes de degré inférieur à $n$. Ainsi, à la limite, l'algèbre de Clifford $C(T_{x_0}^*)$ se comporte comme $\Lambda^\bullet(T_{x_0}^*)$ et on retrouve des formes différentielles. Il reste maintenant à voir quelle forme volume on obtient exactement à la fin. \\

Maintenant, soit $L$ l'opérateur différentiel $D^2$ écrit dans la trivialisation utilisée. Sous le réechelonnement, $r(u,t,v)$ est solution de l'équation de la chaleur 
$$ (\partial_t + u \delta_u L \delta_u^{-1})r(u,t,v) = 0 $$
Cependant, on montre que $u \delta_u L \delta_u^{-1} = K + O(u^{1/2})$, avec $K$ oscillateur harmonique, dont nous connaissons explicitement la solution de l'équation de la chaleur associée, donnée par la formule de Mehler. En utilisant l'unicité de la solution formelle, que l'on obtient après un développement asymptotique de $r(u,t,v)$, on en déduit la formule de Patodi-Gilkey: 
$$ \underset{t \to 0}{\lim} \, \Str(p_t(x,x)) =  \frac{1}{(2i\pi)^{n/2}} \widehat{A}(M) \ch(E/ S)_{[n]} $$
qui permet de conclure. L'objet de cette note est d'appliquer le calcul de Getzler pour donner une légère généralisation du théorème précédent en cohomologie cyclique entière.

\section{Calcul du cocycle JLO pour l'opérateur de Dirac classique}

Soit $M$ une variété riemannienne fermée de dimension paire, $E$ un module de Clifford au-dessus de $M$ muni d'une connexion $\nabla^E$ et $D$ l'opérateur de Dirac associé. Nous calculons maintenant le cocycle JLO du module de Fredholm $\theta$-sommable  $(C^1(M), \Gamma_{L^2}(E), D)$ en cohomologie cyclique entière via le calcul de Getzler. Soient $a_0, \ldots, a_n \in  C^1(M)$, avec $n$ pair. Rappelons que 
$$\mathrm{JLO}_n^t(a_0, \ldots, a_n) = t^{n/2} \int_{\Delta_n} \Str(a_0 e^{-ts_0 D^2}[D,a_1]e^{-ts_1 D^2} \ldots [D,a_n] e^{-ts_n D^2}) ds_0 \ldots ds_n $$
Nous voulons montrer le fait suivant.
\begin{thm*} 
Dans ce cadre, le cocycle $\mathrm{JLO}^t$ restreint à la sous-algèbre dense $C^\infty(M) \subset C^1(M)$ a une limite lorsque $t \to 0$, donnée par la formule suivante :
$$ \mathrm{JLO}_n^t(a_0, \ldots, a_n) \underset{t \to 0}{\longrightarrow} \dfrac{1}{n! \, (2 \mathbf{i} \pi)^{n/2}} \int_M a_0 da_1 \ldots da_n \widehat{A}(M) \ch(E/S)_{[\dim(M) - n]} $$
pour $a_0, \ldots, a_n \in C^\infty(M)$
\end{thm*}

\noindent \textbf{Preuve.}  L'idée est de "ramener" tous les termes $e^{-ts_i D^2}$ "vers la droite" afin d'obtenir
$$e^{-t(s_0 + \ldots + s_n) D^2} = e^{-tD^2} $$
que l'on connaît bien via le calcul de Getzler. Pour ce faire, l'astuce est d'utiliser une technique de commutateurs itérés introduite par Connes et Moscovici \cite{CM}. Pour $a \in C^\infty(M)$, notons 
$$ a(t) = e^{-tD^2} a e^{tD^2} \underset{t \to 0}{\sim} \sum_{k \geq 0} \frac{(-t)^k}{k!} \ad^{(k)}(D^2)(a) $$
où $\ad(D^2)(a) = [D^2,a]$. Cette formule se voit simplement par un développement en série de Taylor en $t=0$. 
Nous noterons également $a^{(k)} := \ad^{(k)}(D^2)(a)$. En particulier, 
$$ e^{-tD^2} a  \underset{t \to 0}{\sim} \sum_{k \geq 0} \frac{(-t)^k}{k!} a^{(k)} e^{-tD^2} $$
Pour $t$ suffisamment petit, 
\begin{eqnarray*}
\mathrm{JLO}_n^t(a_0, \ldots, a_n) & = & t^{n/2} \int_{\Delta_n} \Str \left(a_0 \sum_{k_1 \geq 0} \frac{(-ts_0)^{k_1}}{k_1!} [D,a_1]^{(k_1)} e^{-t(s_0 + s_1)D^2} [D,a_2] \ldots [D,a_n] e^{-ts_n D^2} \right) ds \\
& = & \sum_{k_1, \ldots, k_n \geq 0} \int_{\Delta_n} \frac{(-ts_0)^{k_1}}{k_1!} \ldots \frac{(-ts_n)^{k_n}}{k_1!}\Str \left(a_0 [D,a_1]^{(k_1)} \ldots [D,a_n]^{(k_n)} e^{-t(s_0 + \ldots + s_n) D^2} \right) ds \\
& = & \sum_{k_1, \ldots, k_n \geq 0} \Str (a_0 [D,a_1]^{(k_1)} \ldots [D,a_n]^{(k_n)} e^{-tD^2}) \int_{\Delta_n} \frac{(-ts_0)^{k_1}}{k_1!} \ldots \frac{(-ts_n)^{k_n}}{k_1!} ds
\end{eqnarray*}
L'intégrale ci-dessus vaut $ c(k_1, \ldots , k_n)(-t)^{k_1 + \ldots + k_n} $, où on a noté
$$ c(k_1, \ldots , k_n) = \frac{(-1)^{k_1 + \ldots + k_n}}{k_1 ! \ldots k_n ! (k_1 + 1)(k_1 + k_2 + 2) \ldots (k_1 + \ldots k_n + n)}$$

Remarquons que $[D,a_1]^{(k_1)} \ldots [D,a_n]^{(k_n)}$ est un opérateur différentiel d'ordre au plus $k_1 + \ldots + k_n$. En effet, observons comment agit l'opérateur différentiel $[D^2, [D,f]] = [D^2, c(df)]$ sur les sections de $E$, où $f \in C^{\infty}(M)$ et $c$ est l'action de Clifford. Par la formule de Lichnerowicz, on a une décomposition:
$$ D^2 = \Delta^E + F^{E/S} + \dfrac{1}{4}r_M $$
où $F^{E/S}$ et $r_M$ sont respectivement la courbure tordue et la courbure scalaire. Le laplacien $D^2$ est d'ordre 2, puisque si $(g^{ij})$ est la matrice de la métrique sur $T^*M$, le Laplacien généralisé $\Delta^E$ s'écrit localement sous la forme 
$$ -\sum_{ij} g^{ij} \partial_i \partial_j + \mbox{ termes d'ordre } < 2 $$
de sorte que son symbole principal est la norme sur $T^*M$. Il est alors facile de voir que 
$$ [g^{ij} \partial_i \partial_j , c(da) ] = g^{ij}(\partial_i \partial_j c(da) + \partial_j (c(da)) \partial_i + \partial_i (c(da)) \partial_j) $$
qui est d'ordre au plus $1$, montrant que $[D,a_1]^{(k_1)} \ldots [D,a_n]^{(k_n)}$ est un opérateur différentiel d'ordre au plus $k_1 + \ldots + k_n$. \\

Notant $p_t \in \Gamma(E \boxtimes E^*)$ le noyau de la chaleur de $D^2$, on sait que 

\begin{flushleft}
$t^{n/2}t^{k_1 + \ldots + k_n} \Str (a_0 [D,a_1]^{(k_1)} \ldots [D,a_n]^{(k_n)} e^{-tD^2}) $
\end{flushleft}

\begin{flushright}
 $= \displaystyle{\int_{x \in M} t^{n/2} t^{k_1 + \ldots + k_n} \Str((a_0 [D,a_1]^{(k_1)} \ldots [D,a_n]^{(k_n)})(x) p_t(x,x)) \vert dx \vert}$
\end{flushright}
rappelant que $p_t$ restreint à la diagonale s'identifie à une section de $\End(E)$. Pour effectuer ce calcul, nous allons donc utiliser le calcul de Getzler, en tenant compte des facteurs supplémentaires impliqués. \\

Pour $x_0 \in M$ fixé, on se place dans un système de coordonnées géodésiques centré en $x_0$. Dans un voisinage suffisamment petit, les sections de $E$ s'identifient à des fonctions $C^{\infty}(U, E_{x_0})$ via transport parallèle; si $x, y \in M$ sont tels que $x = \exp_y(v)$ pour un certain $v \in T_yM$, on note $\tau(x,y) : E_y \to E_x$ le transport parallèle le long de la géodésique $t \in [0,1] \longmapsto \exp_y(tv) \in M$ reliant $y$ à $x$, et on pose %Cette fois, la fonction qui va jouer le rôle de $k(t,v)$, où $t > 0$ et $v \in U$, n'est pas $\tau(x_0,\exp_{x_0}(v))p_t(x,\exp_{x_0}(v))$, à la place, on va prendre
\begin{align*}
k(t,v) & = \tau(x_0,\exp_{x_0}(v)) \, t^{n/2}t^{k_1 + \ldots + k_n}(a_0 [D,a_1]^{(k_1)} \ldots [D,a_n]^{(k_n)})(\exp_{x_0}(v)) p_t(\exp_{x_0}(v),x_0) \\
& = \tau(x_0,\exp_{x_0}(v)) \, t^{n/2}t^{k_1 + \ldots + k_n}(a_0 c(da_1)^{(k_1)} \ldots c(da_n)^{(k_n)})(\exp_{x_0}(v)) p_t(\exp_{x_0}(v),x_0) 
\end{align*}
de sorte à voir, grâce au symbole, $k(t,v)$ comme un élément de $C^\infty(U, \End(E_{x_0})) = C^\infty(U, C(V^*) \otimes W) = C^\infty(U, \Lambda^\bullet V^* \otimes W)$, où $V = T_{x_0}M$ et $W = \mathrm{Hom}_{C(V^*)}(S, E_{x_0})$, $S$ étant le module des spineurs associé à l'algèbre de Clifford de $V^*$. Maintenant, rééchelonnons $k(t,v) \in C^\infty(U, \Lambda^\bullet V^* \otimes W)$, i.e on pose pour $u \in (0,1)$
$$ r(u,t,v) = u^{\dim(M)/2} (\delta_u k)(t,v) $$
où 
\[ (\delta_u \alpha)(t,v) = \sum_{i=0}^n u^{-i/2} \alpha(ut,u^{1/2}v)_{[i]} \quad ; \quad \forall \alpha \in C^{\infty}(\mathbb{R}_+ \times U, \Lambda^{\bullet}V^* \otimes \End(W))\]
Rappelons, très important, que ce rééchelonnement est fait de telle sorte à ce qu'on aie 
$$ \boxed{\underset{u \to 0}{\lim} \, r(u,t=1,v=0) = \underset{u \to 0}{\lim} \, \sum_{i=0}^{\dim(M)/2} u^{(\dim(M)/2-i)/2} (a_0(x_0) c(da_1)^{(k_1)}(x_0) \ldots c(da_n)^{(k_n)}(x_0) k_u(x_0,x_0))_{[i]} }$$

Nous connaissons déjà bien le comportement de la partie concernant le noyau de la chaleur lorsque $u \to 0$, il nous faut maintenant examiner le comportement des facteurs 
$$[D,a_i]^{(k_i)} = c(da_i)^{(k_i)} = [D^2, [ D^2 , \ldots[D^2, c(da_i)] \ldots ]]]$$ 
sous l'effet du rééchelonnement. Rappelons que $c(da_i) = \varepsilon(da_i) - \iota(da_i)$. On va regarder le rééchelonnement de $[t D^2, t^{1/2} c(da_i)]$, les facteurs en $t$ proviennent de $t^{n/2}t^{k_1 + \ldots + k_n}$, qu'on incorpore dans les crochets  $[D^2 [ D^2 , \ldots [D^2,c(da_i)] \ldots ]]$. On a
$$ \begin{array}{l}
\delta_u (t^{1/2} c(da_i)) = t^{1/2} (\varepsilon(da_i) - u \iota(da_i)) \\
\delta_u (tD^2)  = ut \delta_u D^2 \delta_u^{-1} = t(K + O(u^{1/2}))
\end{array}$$
où $K = - \sum_i (\partial_i + \frac{1}{4}R_{ij}v_i \partial_j)^2 + F $ est l'oscillateur harmonique. On en déduit alors facilement que 
$$[t D^2, t^{1/2} c(da_i)] = O(u^{1/2})$$
Ainsi, sous le rééchelonnement, $[D,a_i]^{(k_i)}$ s'annule si $k_i \neq 0$ s'annulent lorsque $u \to 0$, par conséquent, d'après la formule de Mehler (cf. par exemple \cite{BGV}), 
$$ \lim_{u \to 0} \, r(u,t,v) = da_1 \wedge \ldots \wedge da_n \wedge (4 \pi t)^{-n/2} \det^{1/2}\left(\frac{tR/2}{\sh(tR/2)}\right) \exp \left(\frac{1}{4t}\left( v, \frac{tR}{2} \coth \left(\frac{tR}{2} \right) v \right) \right) \exp(-tF) $$
Pour $t=1$ et $v=0$, on obtient

\begin{flushleft}
$ \displaystyle{\underset{u \to 0}{\lim} \, \sum_{i=0}^{\dim(M)/2} u^{(\dim(M)/2-i)/2} (a_0(x_0) c(da_1)^{(k_1)}(x_0) \ldots c(da_n)^{(k_n)}(x_0)k_u(x_0,x_0))_{[i]}} $ 
\end{flushleft}

\begin{flushright}
$ = (4 \pi)^{-\dim(M)/2} da_1 \wedge \ldots \wedge da_n \wedge \widehat{A}(M) \exp(-F^{E/S})_{[\dim(M)/2 - n]} $
\end{flushright}
Puis en prenant la supertrace relative, on obtient finalement 
$$ \mathrm{JLO}_n(a_0, \ldots , a_n) = (2i \pi)^{-n/2} \int_M \frac{1}{n!} da_1 \wedge \ldots \wedge da_n \wedge \widehat{A}(M) \ch(E/S)_{[\dim(M) - n]} $$
le terme $1/n!$ étant $c(0, \ldots, 0)$. 
\hfill{$\square$} \\

Le couplage cohomologie cyclique entière/K-théorie donne immédiatement le résultat suivant:

\begin{thm*}
Soit $[e] \in K_0(C^1(M)) \simeq K_0(C^\infty(M))$ représenté par l'idempotent $e \in M_k(C^1(M))$. Alors on a la formule suivante
$$ \Ind(e D_k e) = \frac{1}{(2i\pi)^{\dim(M)/2}}\int_M \ch_{dR}(e)\widehat{A}(M) \ch(E/S) $$
où $\ch_{dR}(e) = \Tr(\exp(e \, de \, de)) $ et $D_k = \mathrm{diag}_k(D, \ldots , D)$.
\end{thm*}

En particulier, pour $e = \textbf{1}$, on retrouve le théorème de l'indice d'Atiyah-Singer.

\end{document}